\documentclass[reqno]{pspum-l}

\usepackage{general,strings}

\numberwithin{equation}{section} \numberwithin{theorem}{section}

\newcommand{\linenopar}{} 

\begin{document}


\title{Tube formulas for self-similar fractals}

\author[M. L. Lapidus]{Michel L. Lapidus}
  \address{University of California, Riverside \\ Department of Mathematics \\ 900 University Ave \\ Riverside, CA 92521}
  \email{lapidus@math.ucr.edu}

\author[E. P. J. Pearse]{Erin P. J. Pearse}
  \address{University of Iowa \\ Department of Mathematics \\ 25L MacLean Hall \\ Iowa City, IA 52246}
  \email{erin-pearse@uiowa.edu}

\thanks{}

\date{\today}



\begin{abstract}
  Tube formulas (by which we mean an explicit formula for the volume of an $\varepsilon$-neighbourhood of a subset of a suitable metric space) have been used in many situations to study properties of the subset. For smooth submanifolds of Euclidean space, this includes Weyl's celebrated results on spectral asymptotics, and the subsequent relation between curvature and spectrum. Additionally, a tube formula contains information about the dimension and measurability of rough sets. In convex geometry, the tube formula of a convex subset of Euclidean space allows for the definition of certain curvature measures. These measures describe the curvature of sets which are not too irregular to support derivatives. In this survey paper, we describe some recent advances in the development of tube formulas for self-similar fractals, and their applications and connections to the other topics mentioned here.
\end{abstract}

\maketitle

\section{Introduction}
\label{sec:Introduction}

This survey article describes advances on the computation of tube formulas for fractal subsets of \bR and of \bRd, and relations to classical results. In particular, we show how the theory of complex dimensions can be used to calculate explicit tube formulas for a large class of self-similar fractals. We also discuss generalizations and connections to forthcoming work.

We begin by discussing results for fractal strings, that is, fractal subsets of \bR. Also, we give the generalization to tube formulas of measures, and show how these results may be applied to the investigation of dimension and measurability of rough sets. In particular, we use the central notion of ``complex dimensions'', a concept which extends the real-valued notion of Minkowski dimension (and Hausdorff dimension, in some cases).
\S\ref{sec:koch} describes the explicit computation of a tube formula for the von Koch snowflake curve, a classical subset of \bRt with fractal boundary. This computation is mainly a ``brute-force'' calculation using elementary geometry coupled with some subtle exploitation of distributional methods.
\S\ref{sec:tilings} contains the construction of a self-similar tiling via an iterated function system, and shows how the tiling enables one to extend the results for fractal strings to self-similar subsets of \bRd.
\S\ref{sec:convex} gives a brief description of some classical results from convex geometry and shows how these relate to the results for self-similar tilings.
Finally, we conclude with a pair of contrasting examples that illustrate the key ideas discussed.

\subsection{Basic concepts}
\label{ssec:basics}

We now present some of the ideas that are used throughout, and introduce basic notation.
\begin{defn}
  \label{def:V of a set}
  Given $\ge > 0$, the \emph{inner \ge-neighbourhood} of a set $A \ci \bRd$,
  $d \geq 1$, is
  \linenopar
  \begin{equation}\label{eqn:def:A_ge}
    A_\ge := \{x \in A \suth dist(x,\del A) \leq \ge\},
  \end{equation}
  where $\del A$ is the boundary of $A$. For a given $A$, we are primarily interested in finding a \emph{tube formula} for $A$, that is, an explicit expression for the $d$-\dimnl Lebesgue measure of $A_\ge$, denoted
  \linenopar
  \[V_{A}(\ge) := \vol[d](A_\ge).\]
\end{defn}

\begin{defn} \label{def:self-similar_system}
  A \emph{self-similar system} is a family $\{\simt_\j\}_{\j=1}^J$ (with $J \geq 2$) of contraction similitudes
  \linenopar
  \[\simt_\j(x) := \simrat_\j A_\j x + a_\j, \q \j=1,\dots,J.\]
  For $\j=1,\dots,J$, we have $0 < \simrat_\j < 1, a_\j \in \bRd$, and $A_\j \in O(d)$, the orthogonal group of rigid rotations in $d$-\dimnl Euclidean space \bRd. The number $\simrat_\j$ is the \emph{scaling ratio} of $\simt_j$. For convenience, assume that
  \linenopar
  \begin{equation}
    \label{eqn:scaling ratio ordering}
    1> \simrat_1 \geq \simrat_2 \geq \dots \geq \simrat_J >0.
  \end{equation}
  When $d=1$, one has only $A_\j = \pm 1$.
\end{defn}

It is well known that there is a unique nonempty compact subset
$\attr \ci \bRd$ satisfying the fixed-point \eqn
  \linenopar
  \begin{equation}
    \label{eqn:fixed-pt_eqn}
    \attr = \simt(\attr) := \bigcup_{\j=1}^J \simt_\j(\attr).
  \end{equation}
This (self-similar) set $F$ is called the \emph{attractor} of \simt. Given a word $w=w_1\dots w_k \in \Wds_k := \{1,2,\dots,J\}^k$, we denote the composition of several similarity mappings by $\simt_w := \simt_{w_k} \comp \dots \comp \simt_{w_2} \comp \simt_{w_1}$.

\section{Fractal strings}
\label{sec:strings}

The essential strategy of fractal strings is to study fractal subsets of \bR by studying their complements.

\begin{defn}
  \label{defn:fractal-string}
  A \emph{fractal string} is any bounded open subset $L \ci \bR$, that is, a countable \seq of open intervals
  \linenopar
  \begin{align} \label{eqn:defn:fractal-string}
    L :=& \{L_n\}_{n=1}^\iy.
  \end{align}
\end{defn}
Although it is not part of the definition, the idea is that the boundary $\del L$ is a set one wants to study. When analyzing $L$ in terms of characteristics which are invariant with respect to rigid motions, the position of an interval $L_n$ within \bR is immaterial, and all pertinent data is stored in the sequence of lengths $\ell_n$ of $L_n$. Indeed, this is the formulation in which the concept is used in \cite{FGCD}, the primary reference for fractal strings.\footnote{This notion was introduced in \cite{LaPo1} and used extensively in \cite{FGCD}, as well as \cite{HeLa}, \cite{La1}, \cite{La2}, \cite{LaMa}, and elsewhere.} For agreement with the sequel, however, we use the inradius.

\begin{defn}
  \label{def:inradius}
  The \emph{inradius} \gr of a set $A$ is
  \linenopar
  \begin{equation}
    \label{eqn:def:inradius}
    \gr = \gr(A)
      := \sup \{\ge > 0 \suth \exists x \text{ with } B(x,\ge) \ci A\}.
  \end{equation}
\end{defn}
It is clear that if $A$ is a bounded set, $A \ci A_\ge$ for sufficiently large \ge. Alternatively, it is apparent that for a fixed $\ge>0$, any sufficiently small set will be entirely contained within its \epsnbd. The notion of \emph{inradius} allows us to see when this phenomenon occurs; its relevance should be clear from Definition~\ref{def:V of a set}.

Obviously, the inradius of an interval is $\gr_n = \frac12 \ell_n = \frac12 \vol[1](L_n)$, i.e., half the length of the interval $L_n \ci \bR$. We then divide by the first (largest) length to normalize the inradii. The result is a \seq of scaling ratios (or scales)  $\{\strrat_n\}$, and if $\genir$ is the length of the largest interval, then $\ell_n = 2\genir\strrat_n$. Use of the terminology ``scaling ratio'' here corresponds to the implicit idea that each interval $L_n \in L$ is congruent to a copy of the largest interval which has been scaled by $\strrat_n$. For a self-similar string, this mapping is explicitly given by $\simt_w$; see \eqref{eqn:scaling-ratio-of-self-sim}. Thus, by the remarks following Definition~\ref{defn:fractal-string}, all pertinent data is stored in the sequence of scaling ratios
  \linenopar
\begin{align} \label{eqn:defn:strings-as-seq-of-inradii}
  \sL :=& \{\strrat_n\}_{n=1}^\iy.
\end{align}
We may take the \seq in nonincreasing order so that $\strrat_1 \geq \strrat_2 \geq \dots > 0$. To avoid trivialities, assume that $\strrat_n >0$ for all $n \in \bN$. Of course, since $L$ is \bdd we have $\sum_{n=1}^\iy \strrat_n < \iy$.

To define a \emph{self-similar string}, consider the self-similar system with similarity mappings $\simt_\j(x) = \simrat_\j A_\j x + a_\j$. The set of scaling ratios of a self-similar string will consists of the collection of all products of scaling ratios of the maps. In particular, if $w$ is a finite word on $\{1,2,\dots,J\}$, then the string will include
\begin{align}
  \label{eqn:scaling-ratio-of-self-sim}
  \strrat_w = \simrat_{w_1} \simrat_{w_2} \dots \simrat_{w_k},
\end{align}
the scaling ratio of $\simt_w = \simt_{w_k} \comp \dots \comp \simt_{w_2} \comp \simt_{w_1}$.
Thus, a self-similar string contains every number that arises as the scaling ratio of a composition of the similarity transformations $\simt_j$. (The first one is always $\strrat_1 = 1$, corresponding to a composition of 0 similarities.) The motivation for this definition is that it corresponds to the set obtained by taking the set-theoretic difference of a self-similar subset of \bR from the smallest closed interval containing it. Alternatively, the set may be constructed by selecting any \bdd open interval $I \ci \bR$ and examining the lengths of the intervals $\{\simt_w(I)\}$, where $w$ runs over all finite words. An example of a self-similar string is given by the Cantor string in Example~\ref{exm:Cantor string} just below.

A \emph{generalized fractal string} is a locally finite Borel measure on $(\ge,\iy)$, where $\ge>0$, and is denoted $\gh=\gh_\sL$. Such a string may not have a geometric realization. The motivation for this generalization lies in the flexibility of working in the measure-theoretic framework, and in certain specific applications. In this context, an ordinary fractal string as given in the previous definitions (self-similar or not) corresponds to a sum of Dirac masses $\gd_x$, each located at a reciprocal of one of the numbers $\strrat_n$. That is, an ordinary fractal string may be written
  \linenopar
  \begin{equation}
    \label{eqn:ghL-as-dirac-sum}
    \ghL = \sum_{n=1}^\iy \gd_{1/\strrat_n}.
  \end{equation}

\begin{defn}\label{defn:scaling-zeta-fn}
  The \emph{scaling zeta function} \gzs of a fractal string is the Mellin transform of the measure $\ghL$:
  \linenopar
  \begin{align}\label{eqn:defn:scaling-zeta-fn}
    \gzs(s) := \int_0^\iy x^{-s} \,d\ghL(x), \qq s \in \bC.
  \end{align}
  In \cite{FGCD}, this is called the geometric zeta \fn of \sL. We have chosen the current terminology to agree with the latter sections of this paper and instead say that the \emph{geometric zeta function} of \gzL is given by
  \linenopar
  \begin{align}\label{eqn:defn:geometric-zeta-fn-string}
    \gzL(\ge,s) := \gzs(s)\frac{(2\ge)^{1-s}}{s(1-s)}.
  \end{align}
\end{defn}
The \fn \gzL factors into \gzs and a term which contains geometric data about open intervals, although this will not become clear until the discussion of the higher-dimensional tilings in \S\ref{sec:tilings}. In the case when \sL is self-similar, the function $\gzL:\gW \to \bC$ may be meromorphically continued to all of \bC. Otherwise, it is well defined in some half-plane $\{\Re s > a\}$ and may be continued analytically. To demarcate the domain of \gzL formally, we introduce the following definition.

\begin{defn}
  \label{def:screen}
  Let $f:\bR \to \bR$ be a \bdd Lipschitz \cn \fn. Then the \emph{screen} is $S = \{f(t) + \ii t \suth t \in \bR\}$, the graph of $f$ with the axes interchanged. The region to the right of the screen is the \emph{window} $W := \{z \in \bC \suth \Re z \geq f(\Im z)\}$. We choose $f$ so that $S$ does not pass through any poles of \gzL, and \gzL has a \mero extension to some \nbd of $W$.
\end{defn}

\begin{defn}
  \label{def:complex-dimensions}
   The \emph{(complex) scaling dimensions} of \sL are poles of \gzL. The poles which lie in the window are called the \emph{visible scaling dimensions}:
  \linenopar
  \begin{align}
    \label{eqn:def:visible-scaling-dimns}
    \D_\sL(W) = \{\gw \in W \suth \lim_{s \to \gw} |\gzL(s)| = \iy\}.
  \end{align}
\end{defn}

The screen and window are useful for many purposes. In addition to demarcating the domain of \gzL, they allow one to make precise statements about the growth of \gzL, and they allow for certain approximation arguments that make the proof of the next theorem possible. In particular, one has the following definition.

\begin{defn}\label{def:languid}
  One says that \gzL (or just \sL) is \emph{languid} if it \sats certain mild growth \conds on $S$, and along a \seq of \horz lines in $W$. For the precise nature of these \conds, see \cite[Def.~5.2]{FGCD}.
\end{defn}

\begin{theorem}[Tube formula for fractal strings {\cite{FGCD}, Thm.~8.1}]
  \label{thm:tube_formula_for_fractal_strings}
  Let $\ghL$ be a fractal string with geometric zeta \fn \gzL and assume that \gzL is languid. Then we have a tube formula
  \linenopar
  \begin{equation}
    \label{eqn:1-dim_tube_formula}
    V_\sL(\ge)
    = \negsp[15] \sum_{\gw \in \D_\sL(W)} \negsp[15] \res[s=\gw]{\gzL(\ge,s)}
     + \{2\ge \gzh(0)\} + \sR(\ge).
  \end{equation}
  Here the term in braces is only included if $0 \in W \less \D_\gh(W)$, and
  the error term is
  \linenopar
  \begin{equation}
    \label{eqn:1-dim_error_term}
    \sR(\ge)
    = \frac1{2\gp \ii} \int_S \gzL(\ge,s) \, ds
    = O(\ge^{1-\sup S}), \qq \text{as } \ge \to 0^+.
  \end{equation}
  When the string is also self-similar, one may take $W=\bC$ and $\sR(\ge) \equiv 0$.
\end{theorem}
If we denote the poles of the scaling zeta \fn separately by \Ds, then the result above may be rewritten
  \linenopar
  \begin{equation}
    \label{eqn:1-dim_tube_formula-conceptual}
    V_\sL(\ge)
    = \negsp[20] \sum_{\gw \in \Ds(W) \cup \{0\}} \negsp[20] \res[s=\gw]{\gzL(\ge,s)} + \sR(\ge).
  \end{equation}
Formula \eqref{eqn:1-dim_tube_formula} was originally obtained as a \dist acting on smooth \fns with compact support in $(0,\iy)$. However, it has since been obtained in a pointwise fashion in \cite[Thm.~8.7]{FGCD} under only slightly more restrictive \conds (which are always \satd in the case of self-similar strings), so this technicality need not be emphasized.

\begin{exm}\label{exm:Cantor string}
  The complement of the usual Cantor set in the unit interval $[0,1]$ consists of open intervals with lengths $\{\frac13, \frac19, \frac19, \frac1{27}, \dots,\}$. This self-similar fractal string has a largest interval of length $\frac13$ and scaling ratios $3^{-k}$ appearing with multiplicity $2^k$. Consequently, the Cantor string $\CS$ may be written
  \linenopar
  \begin{align}\label{eqn:cantor-measure}
    \gh_{\CS} = \sum_{k=0}^\iy 2^k \gd_{3^k},
  \end{align}
  and its scaling zeta \fn is
  \linenopar
  \begin{align}\label{eqn:cantor-zetas}
    \gzs(s) = \int_0^\iy x^{-s} \,d\gh_{\CS}(x)
    = \sum_{k=0}^\iy \left(\frac{2}{3^s}\right)^k = \frac{1}{1-2\cdot3^{-s}},
  \end{align}
  while its tube formula is
  \linenopar
  \begin{align}\label{eqn:cantor-tube-formula}
    V_\sC(\ge)
    &= \frac1{3\log3} \sum_{n \in \bZ} \left(\frac{1}{D+\ii n\per}-\frac1{D-1+\ii n\per}\right)\left(\frac{\ge}{\genir}\right)^{1-D-\ii n\per} - 2\ge,
  \end{align}
  with $\genir=\tfrac16$ is the inradius of the largest interval, $D=\log_32$ is the Minkowski dimension, and $\per=2\gp/\log3$ is a constant called the \emph{oscillatory period}.
\end{exm}

\begin{defn}\label{defn:Minkowski-dimn}
  The \emph{Minkowski dimension} of $A \ci \bRd$ is \linenopar
  \begin{align}\label{eqn:defn:Minkowski-dimn}
    D = \dim_M A := \inf\{t \geq 0 \suth V_A(\ge) = O(\ge^{d-t}), \text{ as } \ge \to 0^+\},
  \end{align}
  and is also frequently called the box dimension. For a string, $\dim \sL := \dim \del L$.
\end{defn}

The complex dimensions can be thought of as a generalization of Minkowski dimension because of the following result of \cite{La2}, which also appears as \cite[Thm.~1.10]{FGCD}.
\begin{theorem}
  For a string \sL, the scaling zeta \fn \gzs converges on the half-plane $\{s \in \bC \suth \Re s > \gs\}$ if and only if $\gs \geq \dim \sL$.
\end{theorem}

\begin{remark}\label{rem:self-similar_zeta}
  In formula \eqref{eqn:cantor-zetas}, one discovers the pleasant surprise that \gzs has a nice closed form. In fact, it is shown in \cite[Thm.~2.4]{FGCD} that all self-similar strings have a scaling zeta \fn of the form
  \linenopar
  \begin{align}\label{eqn:self-similar_zeta}
    \gzs(s) = \frac1{1-\sum_{j=1}^J \simrat_j^s}.
  \end{align}
  This remark remains true for the self-similar tilings discussed in \S\ref{sec:tilings}. The number-theoretic and measure-theoretic implications of this result are extensive; they are studied at length in \cite[Ch.~2--3]{FGCD}.
\end{remark}

\section{The von Koch curve}
\label{sec:koch}

\begin{figure}
  \includegraphics{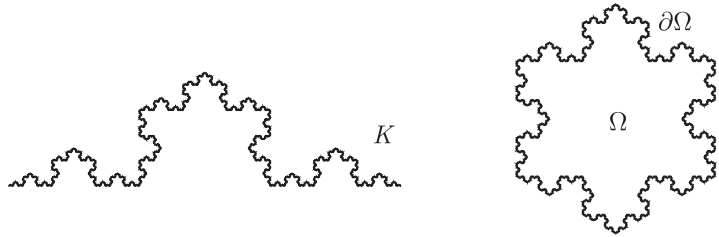}
  \caption{The Koch curve $K$ and Koch snowflake domain \gW.}
  \label{fig:koch-snowflake} \centering
\end{figure}

The Koch curve $K$ is the attractor for the self-similar system
  \linenopar
\begin{align}
  \label{eqn:koch-system}
  \simt_1(z) := \gr \cj{z} \q \text{and} \q
  \simt_2(z) := (1-\gr)(\cj{z}-1)+1,
\end{align}
where $\gr = \tfrac{1}{2} + \tfrac{1}{2\sqrt{3}}i$ and $\ii=\sqrt{-1}$. Consequently, $K$ is a self-similar fractal with Minkowski dimension $D := \log_3 4$; see Figure~\ref{fig:koch-snowflake}. In this section, we describe how the tube formula for the Koch snowflake was calculated directly.

\begin{figure}
    \includegraphics{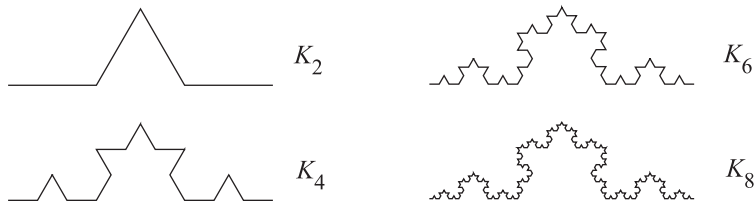}
    \caption{Four early stages in the geometric construction of $K$.}
    \label{fig:koch-pre-curves} \centering
  \end{figure}

\begin{figure}
    \includegraphics{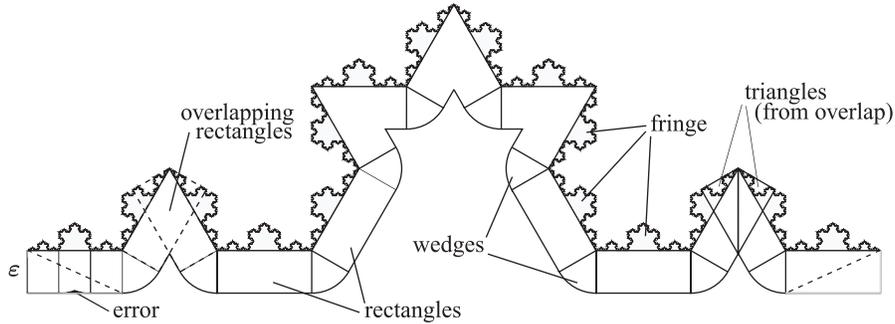}
    \caption{An approximation to the inner \ge-\nbd of the Koch
    curve, based on the curve $K_2$ from Figure~\ref{fig:koch-pre-curves}. }
    \label{fig:midneighborhood} \centering
  \end{figure}

\begin{figure}
    \includegraphics{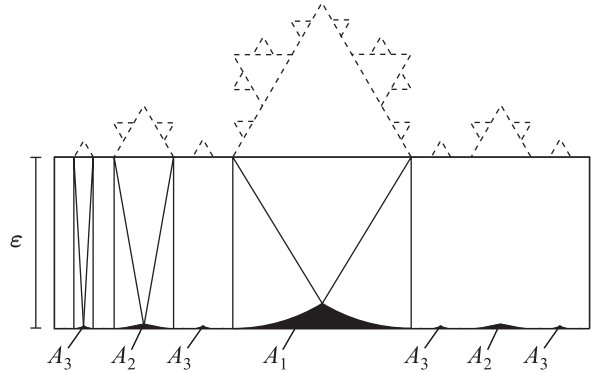}
    \caption{An error block for $K_n$. The central third of
    the block contains one large isosceles triangle, two wedges,
    and the trianglet $A_1$. Figure~\ref{fig:forming-blocks} contains 4 complete copies of this figure, and 12 partial copies of it.}
    \label{fig:error-block} \centering
  \end{figure}

In \cite{KTF}, the authors compute a tube formula for the Koch snowflake \gW by considering a \seq of curves $K_n \to K$, where convergence is with respect to Hausdorff metric. Some representative terms of this \seq are illustrated in Figure~\ref{fig:koch-pre-curves}.
The tube formula is obtained as a limit of tube formulas $V_{K_n}$ obtained by approximation for each $K_n$, as illustrated for $K_2$ in Figure~\ref{fig:midneighborhood}. Three copies of this figure are fitted together to form the inner \nbd of the snowflake, as is indicated by the dashed lines at either end.
This region is only an approximation of course, as one side of each rectangle should be replaced with a fractal curve. Fig.~\ref{fig:error-block} shows how this error is incurred and how it inherits a Cantoresque structure from the triadic character of the Koch curve.
Let us refer to the dark region in Fig.~\ref{fig:error-block} as an \emph{error block} and each connected component $A_k$ of it as a \emph{trianglet}. Without taking the error blocks into account, the \ge-\nbd of the Koch curve has approximate area
\linenopar
\begin{equation}\label{fn:preliminary-area}
  \preV(\ge) = \ge^{2-D} 4^{-\{x\}}
    \left(\tfrac{3\sqrt{3}}{40}9^{\{x\}} + \tfrac{\sqrt{3}}{2} 3^{\{x\}}
    + \tfrac16\left(\tfrac{\gp}{3} - \sqrt{3}\right)\right)
    - \tfrac{\ge^2}3 \left(\tfrac{\gp}{3} + 2\sqrt{3}\right).
\end{equation}
Here, $x := -\log_3(\ge\sqrt{3})$, and we use $x=[x] + \{x\}$ to denote the decomposition of $x$ into its integer and fractional parts. That is, $[x]$ is an integer and $0 \leq \{x\} < 1$. Furthermore, $x$ is related to $n$ (the index of $K_n$) by $n = n(\ge) = [\log_3 \frac{1}{\ge\sqrt{3}}] = [x]$. Thus, the level $n$ of the approximation is determined by \ge, with $n \to \iy$ as $\ge \to 0$.

\begin{figure}
    \includegraphics{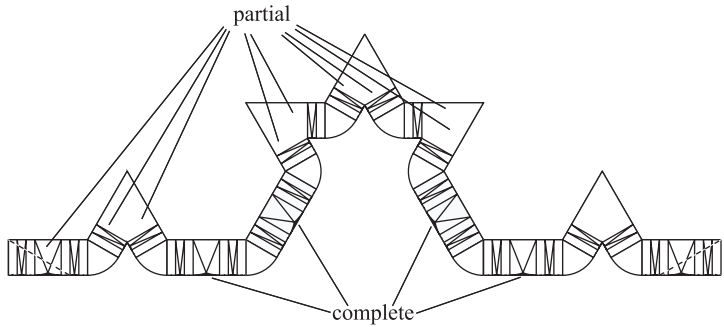}
    \caption{Error block formation. The ends are counted as partial because three copies of this illustration will be added together to make the entire snowflake.}
    \label{fig:forming-blocks}
  \end{figure}

The tube formula for \gW will be obtained by summing the areas $A_k$ of the trianglets, and multiplying by the number of error blocks occurring for a given approximation, i.e., for a particular value of \ge. Unfortunately, this number is not easy to express, due to the existence of partial error blocks. See Figure~\ref{fig:forming-blocks} for a visual explanation of what is meant by \emph{partial error blocks} and \emph{complete error blocks}. 
The number of complete error blocks can be readily counted with a simple formula, but the portion of a partial error block that exists for a given value of \ge is rather ornery, and so we denote this quantity by $h(\ge)$. We do not know $h(\ge)$ explicitly, but we do know by the self-similarity of $K$ that it has multiplicative period 3; i.e., $h(\ge) = h(\tfrac\ge3)$. It is a theorem of \cite{LlWi:fractal-euler-char} that for a self-similar set, the \ge-\nbd must be either: (i) a finite union of convex sets for every value of \ge, or (ii) not a finite union of convex sets for any value of \ge. We conjecture that for every self-similar set of the latter type, the tube formula involves a multiplicatively periodic \fn analogous to our $h(\ge)$.
  \begin{figure}
    \includegraphics{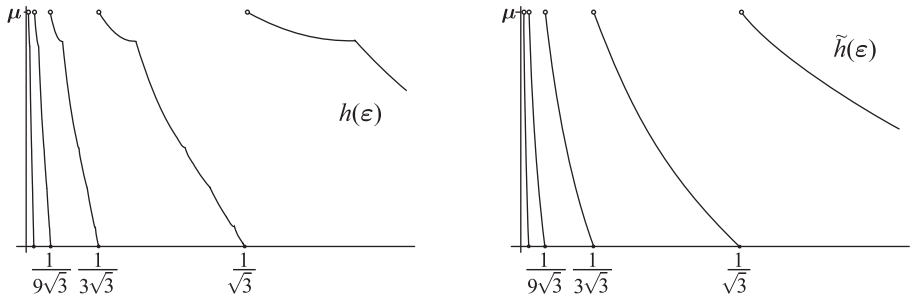}
    \caption{A comparison between the graph of the Cantor-like \fn $h$ and the graph of an approximation $\htil = \gm \cdot (-[x]-x)$, where \gm is a constant and $x=-\log_3 (\ge \sqrt 3)$.}
    \label{fig:h_vs_h} \centering
  \end{figure}
Once the error has been taken into account, one obtains \cite[Thm.~5.1]{KTF} as follows.

\begin{theorem}
  \label{thm:preview-of-final-result}
  The area of the inner \ge-\nbd of the Koch snowflake is given by the following tube formula:
  \linenopar
  \begin{equation}
    \label{eqn:preview-of-corollary}
    V(\ge) = \sum_{n \in \bZ} \gf_n \ge^{2-D-in\per} + \sum_{n \in \bZ} \gy_n \ge^{2-in\per},
  \end{equation}
  for suitable constants $\gf_n, \gy_n$ which depend only on $n$ and are expressed in terms of the Fourier \coeffs $g_\ga$ of $h(\ge)$.
\end{theorem}
To see the full form of \eqref{eqn:preview-of-corollary}, please see Remark~\ref{rem:comparison of Koch tiling tube to Koch tube} at the end of this paper. The result \eqref{eqn:preview-of-corollary} was obtained at a time when the theory of complex dimensions was entirely restricted to fractal strings as outlined in \S\ref{sec:strings} (except for the conjectures expressed in \cite[Ch.~12]{FGCD}). Part of the motivation for proving Theorem~\ref{thm:preview-of-final-result} was to get an idea of what the complex dimensions of \gW might look like. Reasoning by analogy, one would deduce from \eqref{eqn:preview-of-corollary} that the complex dimensions of \gW are obtained from the exponents appearing in \eqref{eqn:preview-of-corollary}; in particular, that each exponent is of the form $2-\gw$, where \gw is a complex dimension of \gW. This led to the prediction that the complex dimensions of the Koch snowflake are
  \linenopar
\begin{equation}
  \label{eqn:complex-dim-of-koch}
  \D_{\del\gW} = \{D+in\per \suth n \in \bZ\}
    \cup \{in\per \suth n \in \bZ\}.
\end{equation}
The reader will find in Example~\ref{exm:The Koch Tiling} and Remark~\ref{rem:comparison of Koch tiling tube to Koch tube} that this is not far wrong.

\section{Fractal sprays and self-similar tilings}
\label{sec:tilings}

A fractal spray is the higher-\dimnl counterpart of a fractal string. See \cite[\S~1.4]{FGCD} for a discussion; this idea also appears earlier in \cite{LaPo2}, \cite{La2}, and \cite{La3}.

\begin{defn}\label{def:fractal-spray}
  Let $\gen \ci \bRd$ be a nonempty \bdd open set, which we will call the \emph{generator}, and let $\sL = \{\strrat_n\}$ be a fractal string. Then a \emph{fractal spray} is a bounded open subset of \bRd which is the disjoint union of open sets $\tile_n$ for $n=1,2,\dots$, where each $\tile_n$ is congruent to $\strrat_n \gen$, the homothetic of \gen by $\strrat_n$.
\end{defn}
Thus, any fractal string can be thought of as a fractal spray on the basic shape $\gen=(0,1)$, the unit interval. Every self-similar system is naturally associated to a certain fractal spray called the self-similar tiling.

\begin{defn}\label{def:self-similar-tiling}
  The \emph{self-similar tiling} \tiling corresponding to a self-similar system \simt is a fractal spray where the fractal string and generators are defined as follows. The string is the collection of all finite products of the scaling ratios $\{\simrat_1, \dots, \simrat_J\}$ of the self-similar system \simt. Let $\hull = [\attr]$ be the convex hull of \attr, and let $\inter \hull$ be its interior. Then the generators are the connected open sets $\gen_q$ in the disjoint union
  \linenopar
  \begin{align}\label{eqn:defn:generators}
    \inter \hull \less \simt(\hull) = \gen_1 \cup \dots \cup \gen_Q.
  \end{align}
\end{defn}
When there is more than one generator, it is more accurate to think of the tiling as a union of fractal sprays, one for each generator. Indeed, for purposes of computing the tube formula, it is easiest to deal with each part separately and then obtain the final result by adding the contributions from each generator. For this reason, the tube formulas below are all phrased for a spray or tiling with one generator. In fact, we cannot currently exclude the possibility that there may exist examples for which $Q = \iy$. However, we have been unable to construct such an example, and for the time being we assume $Q < \iy$.

\begin{figure}
  \includegraphics{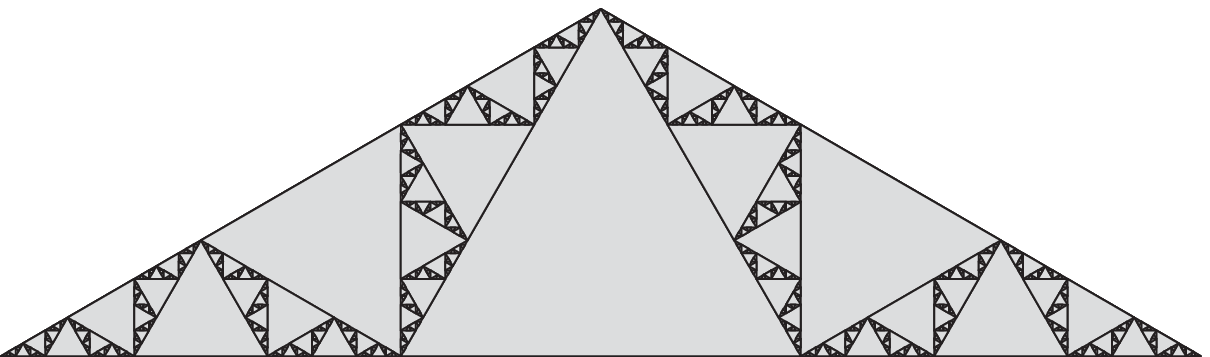}
  \caption{\captionsize The Koch tiling \sK. The generator is the largest equilateral triangle, and the system is $\{\simt_1(z) = \gx \cj{z}, \simt_2(z) = (1-\gx)(\cj{z}-1)+1$, for $z \in \bC$, where $\gx = \tfrac{1}{2} + \tfrac{1}{2\sqrt{3}}\ii$.}
  \label{fig:koch-tiled-in}
  \centering
\end{figure}

The term ``self-similar tiling'' is used here in a sense quite different from the one often encountered in the literature. In particular, the tiles themselves are neither self-similar nor are they all of the same size; in fact, the tiles are typically simple polyhedra. Moreover, the region being tiled is the complement of the self-similar set \attr within its convex hull, rather than all of \bRd; see Figure~\ref{fig:koch-tiled-in}. In fact, it is shown in \cite{SST} that when \simt \sats the tileset condition, the collection $\{\simt_w(\gen_q)\}_{w,q}$ forms an open tiling of $\hull \less \attr$. We now define these terms.

\begin{defn} \label{defn:open-tiling}
  An \emph{open tiling} of $A \in \bRd$ is a collection of nonempty connected open sets $\{A_n\}_{n=1}^\iy$ \st
  \linenopar
  \begin{align}
    \text{(i)}&\; \cj{A} = \cj{\bigcup\nolimits_{n=1}^\iy {A_n}}, \qq \text{and}\qq
    \text{(ii)}\; A_n \cap A_m = \es \text{ for } n \neq m.
  \end{align}
  Figure~\ref{fig:koch-tiled-in} shows a tiling by open sets which are equilateral triangles.
\end{defn}

\begin{defn}
  \label{defn:tileset-condition}
  The system \simt \sats the \emph{tileset condition} iff $\tileset \nsubseteq \simt(\hull)$ and
  \linenopar
  \begin{equation}
    \label{eqn:tileset-condition}
      \inter \simt_j(\hull) \cap \inter \simt_\ell(\hull) = \es,
      \qq j \neq \ell.
  \end{equation}
\end{defn}
The tileset \cond is a separation \cond which is similar to the ``open set \cond'', but stronger. The open set \cond is \satd when there exists a ``feasible open set'' $U$ which has the property that $\simt_j(U) \ci U$, and $\simt_j(U) \cap \simt_\ell(U)$ for all $j,\ell$. To see that the tileset \cond implies this, let $U$ be the the interior of $\hull = [\attr]$. For a counterexample to the converse, see \cite[Example~3.8]{SST}.

For self-similar tilings, and more generally for fractal sprays, the scaling zeta \fn is defined just as it was for fractal strings. However, the geometric zeta \fn becomes much more complicated, due to the multifarious possibilities for the geometry of the generator. The resulting technicalities can be ignored in many cases, however, for example when the generator is sufficiently simple. This motivates the \cond ``Steiner-like'' in the definition just below. A bounded open set in \bRd is said to be Steiner-like if its inner parallel volume $V_\gen(\ge)$ admits a ``polynomial-like'' expansion in \ge of degree at most $d$. More precisely, we have the following definition (see \S\ref{sec:convex} for an explanation of this choice of terminology).
\begin{defn}
  \label{def:Steiner-like}
  A \bdd open set $\gen \ci \bRd$ is \emph{Steiner-like} iff
  for $0 \leq \ge \leq \gr(\gen)$ its inner tube formula may be written
  \linenopar
  \begin{align}\label{eqn:def-prelim-Steiner-like-formula}
    V_\gen(\ge) = \sum_{k=0}^{d} \crv_k(\gen,\ge) \ge^{d-k},
  \end{align}
  where each coefficient \fn $\crv_k(\gen,\cdot)$ is assumed to
  be a bounded and locally integrable function of \ge for which
  \linenopar
  \begin{equation}\label{eqn:tractability-of-kappas}
    \lim_{\ge \to 0^+} \crv_k(\gen,\ge)
  \end{equation}
  exists, and is both positive and finite.
\end{defn}

\begin{defn}
  \label{def:diphase-and-pluriphase}
  In the special case when \gen is a Steiner-like set whose terms $\crv_k$ are constant, we say the set is \emph{diphase}. This terminology refers to the fact that its tube formula is written piecewise with only two cases: $\ge \leq \gr(\gen)$ and $\ge > \gr(\gen)$. If \gen is Steiner-like and the \fns $\crv_k$ are piecewise constant on the interval $[0,\gr(\gen)]$, then we say \gen is \emph{pluriphase}. Thus, diphase is a special case of pluriphase.
\end{defn}

\begin{conj}
  The class of pluriphase sets includes convex sets and polyhedra.
\end{conj}
This result may exist in the literature, but we have been unable to find it. It is simple to show that the class of diphase sets includes spheres and regular polyhedra. In particular, the examples of \S\ref{sec:examples} are both diphase. There are examples of other types of diphase sets, but they resist easy description. Current research is attempting to characterize diphase and pluriphase sets in terms of other geometric properties.

In the case when \gen is pluriphase, the geometric zeta \fn is relatively simple to write down. To avoid obscuring the exposition, we give only the diphase case here (although in view of Definition~\ref{def:visible-complex-dimensions-of-a-fractal-spray} it is worth noting that the factors $\frac1{s-\i}$ also appear in the more general case).

\begin{defn}
  \label{def:geometric-zeta-fn-of-a-fractal-spray}
  The \emph{geometric zeta function} of a fractal spray with a diphase generator is
  \linenopar
  \begin{align}
    \label{eqn:diphase-geometric-zeta}
    \gzT(\ge,s) :=& \ge^{d-s} \gzh(s) \sum_{\i=0}^d \frac{\genir^{s-\i}}{s-\i} \crv_{\i}.
  \end{align}
\end{defn}

\begin{defn}
  \label{def:visible-complex-dimensions-of-a-fractal-spray}
  The set of \emph{visible complex dimensions of a fractal spray} is
  \linenopar
  \begin{align}
    \label{eqn:spray-dimensions}
    \DT(W) := \D_\gh(W) 
              \cup \intDim.
  \end{align}
  Thus, $\DT(W)$ consists of the visible scaling dimensions and the ``integer dimensions'' of \gen, and contains all the poles of \gzT (viewed as a \fn of $s$).
\end{defn}

\begin{theorem}[Tube formula for fractal sprays]
  \label{thm:fractal-spray-tube-formula}
  Let \gh be a fractal spray on the Steiner-like generator \gen, with generating inradius $\genir = \gr(\gen) > 0$. If \gzT is languid, then we have a tube formula
  \linenopar
  \begin{align}
    \label{eqn:main-result}
    V_\tiling(\ge) &= \sum_{\gw \in \DT(W)}
    \res{\gzT(\ge,s)} + \sR(\ge),
  \end{align}
  where the sum ranges over the set $\DT(W)$. Here, the error term $\sR(\ge)$ is
  \linenopar
  \begin{align}
    \label{def:main-error-term}
    \sR(\ge) 
    = \frac1{2\gp \ii}\int_S \gzT(\ge,s) \, ds,
    =  O(\ge^{d-\sup S}), \qq \text{as } \ge \to 0^+.
  \end{align}
\end{theorem}
Theorem~\ref{thm:fractal-spray-tube-formula} was first obtained \distly in \cite{Thesis} and appears in improved form in \cite{TFCD}. Since then, a pointwise proof has been obtained in \cite{Pointwise}. It is important to note that \gzT is languid for all self-similar tilings. As we have not been able to construct a self-similar tiling which \sats the tileset \cond but fails to be Steiner-like, Theorem~\ref{thm:fractal-spray-tube-formula} applies to all known examples. \emph{A fortiori}, it is possible to show that all self-similar tilings automatically \sat a more stringent condition (called \emph{strongly languid} in \cite[Def~5.2]{FGCD}) that enables one to take $W = \bC$ and $\sR(\ge) \equiv 0$; see \cite[Thm.~8.4]{TFCD}. Thusly one obtains the following special case of Theorem~\ref{thm:fractal-spray-tube-formula}.

\begin{cor}
  \label{cor:self-similar,_one_generator}
  Let \tiling be a self-similar tiling with a diphase generator \gen. If $\gzT(s)$ has only simple poles, then
  \linenopar
  \begin{align}
    \label{eqn:single_gen_result}
    V_\tiling(\ge)
    &= \sum_{\gw \in \Ds} \sum_{\i=0}^{d}  \res{\gzs(s)} \ge^{d-\gw}
      \tfrac{\genir^{\gw-\i}}{\gw-\i} \crv_\i
      + \sum_{\i=0}^{d-1} \crv_\i \gzs(\i) \ge^{d-\i}.
  \end{align}
\end{cor}
By comparing \eqref{eqn:single_gen_result} with \eqref{eqn:1-dim_tube_formula}, it is easy to see how the tube formula for fractal sprays extends the results for fractal strings. Additionally, it extends classical results for convex sets, as outlined in the next section; see especially \eqref{eqn:steiner_conceptual}--\eqref{eqn:main_cor_conceptual} and the surrounding discussion. 
\section{Convex geometry and the curvature measures}
\label{sec:convex}

In order to explain the connections with convex geometry, we give a brief encapsulation of Steiner's classical result. Here, we denote the Minkowski sum of two sets in \bRd by
  \linenopar
\begin{align*}
  A+B = \{x \in \bRd \suth x=a+b \text{ for } a \in A, b \in B\}.
\end{align*}

\begin{theorem}\label{thm:Steiners-formula-basic}
  If $B^d$ is the $d$-\dimnl unit ball and $A \ci \bRd$ is convex, then the $d$-\dimnl volume of $A+\ge B^d$ is given by
  \linenopar
  \begin{align}
    \vol[d](A+\ge B^d) = \sum_{\i=0}^{d} \gm_\i(A) \vol[d-\i](B^{d-\i}) \ge^{d-\i},
      \label{def:Steiner formula with inv meas}
  \end{align}
  where $\gm_i$ is the renormalized $i$-\dimnl \emph{intrinsic volume} and $\vol[j](A)$ is the $j$-\dimnl Lebesgue measure.
\end{theorem}

Up to some normalizing constant, the $i$-\dimnl intrinsic volume is the same thing as the \ith total curvature or \nth[(d-i)] Quermassintegral. This valuation $\gm_i$ can be defined via integral geometry as the average measure of orthogonal projections to $(d-i)$-\dimnl subspaces; see Chap.~7 of \cite{Rota}. For now, we note that (up to a constant), there is a correspondence
\begin{center}
\begin{tabular}{rrlrrl}
  $\gm_0$ & $\sim$ &Euler characteristic, \qq \qq & $\gm_{d-1}$ & $\sim$ &surface area, \\
  $\gm_1$ & $\sim$ &mean width, & $\gm_d$ & $\sim$ &volume,
\end{tabular}
\end{center}
see \cite[\S4.2]{Schn2} for more. We have chosen the term ``Steiner-like'' for Def.~\ref{def:Steiner-like} \bc the intrinsic volumes satisfy the following properties:
  \linenopar
\begin{enumerate}[(i)]
  \item each $\gm_\i$ is \homog of degree $\i$, so that for $x >0,$
  \label{itm:homogeneity of mu}
  \begin{align}
    \label{eqn:homogeneity of mu}
    \gm_\i\left(x A\right) = \gm_\i(A) \, x^\i, \text{ and}
  \end{align}
\item each $\gm_\i(A)$ is rigid motion invariant, so that
  \label{itm:translation invariance of mu}
  \begin{align}
    \label{eqn:translation invariance of mu}
    \gm_\i\left(T(A)\right) = \gm_\i(A),
  \end{align}
  for any (affine) isometry $T$ of \bRd.
\end{enumerate}
Note that \eqref{def:Steiner formula with inv meas} gives the volume of the set of points which are within \ge of $A$, including the points of $A$. If we denote the \emph{exterior \epsnbd} of $A$ by
  \linenopar
\begin{align*}
  A_\ge^{ext} := (A+\ge B^d) \less A = \{x \suth d(x,A) \leq \ge, x \notin A\},
\end{align*}
then it is immediately clear that omitting the \nth[d] term gives
  \linenopar
\begin{align}
  \vol[d](A_\ge^{ext}) = \sum_{\i=0}^{d-1} C_\i(A) \ge^{d-\i}
      \label{def:Steiner ext formula with inv meas}
\end{align}
with $C_i(A) = \gm_\i(A) \vol[d-\i](B^{d-\i})$. The intrinsic volumes $\gm_i$ can be localized and understood as the \emph{curvature measures} introduced in \cite{Fed} and described in Ch.~4 of \cite{Schn2}. In this case, for a Borel set $\gb \ci \bRd$, one has
  \linenopar
\begin{align}
  \vol[d]\{x \in A_\ge^{ext} \suth p(x,A) \in \gb\}
  = \sum_{\i=0}^{d-1} C_\i(A,\gb) \ge^{d-\i},
    \label{def:Steiner ext formula localized}
\end{align}
where $p(x,A)$ is the metric projection of $x$ to $A$, that is, the closest point of $A$ to $x$. \emph{In fact, the curvature measures are obtained axiomatically in \cite{Schn2} as the \coeffs of the tube formula, and it is this approach that we hope to emulate in our current work.} In other words, we believe that $\crv_i$ may also be understood as a (total) curvature, in a suitable sense, and we expect that $\crv_i$ can be localized as a curvature measure. A more rigorous formulation of these ideas is currently underway in \cite{FCM}. Caveat: the description of $\crv_\i$ given in the conditions of Def.~\ref{def:Steiner-like} is intended to emphasize the resemblance between $\crv_\i$ and $C_i$. However, $\crv_\i$ may be signed (even when \gen is convex and $\i=d-1,d$) and is more complicated in general. In contrast, the curvature measures $C_\i$ are always positive for convex bodies.

The primary reason we have worked with the inner \ge-\nbd instead of the exterior is that it is more intrinsic to the set; it makes the computation independent of the embedding of \tiling into \bRd. At least, this should be the case, provided the `curvature' terms $\crv_\i$ of Def.~\ref{def:diphase-and-pluriphase} are also intrinsic. As a practical bonus, working with the inner \ge-\nbd allows us to avoid potential issues with the intersections of the \ge-\nbds of different components.

In \cite{Fed}, Federer unified the tube formulas of Steiner (for convex bodies, as described in Ch.~4 of \cite{Schn2}) and of Weyl (for smooth submanifolds, as described in \cite{BeGo}, \cite{Gr} and \cite{We}) and extended these results to sets of \emph{positive reach}.\footnote{A set $A$ has \emph{positive reach} iff there is some $\gd>0$ \st any point $x$ within \gd of $A$ has a unique metric projection to $A$, i.e., that there is a unique point $A$ minimizing $dist(x,A)$. Equivalently, every point $q$ on the boundary of $A$ lies on a sphere of radius \gd which intersects $\del A$ only at $q$.}
It is worth noting that Weyl's tube formula for smooth submanifolds of \bRd is expressed as a polynomial in \ge with \coeffs defined in terms of curvatures (in a classical sense) that are \emph{intrinsic} to the submanifold \cite{We}. See \S6.6--6.9 of \cite{BeGo} and the book \cite{Gr}.
Federer's tube formula has since been extended in various directions by a number of researchers in integral geometry and geometric measure theory, including [\citen{Schn1},\citen{Schn2}], [\citen{Za1},\citen{Za2}], [\citen{Fu1},\citen{Fu2}], \cite{Stacho}, and most recently (and most generally) in \cite{HuLaWe}. The books \cite{Gr} and \cite{Schn2} contain extensive endnotes with further information and many other references.

To emphasize the present analogy, consider that Steiner's formula \eqref{def:Steiner ext formula with inv meas} may be rewritten
  \linenopar
\begin{align}
  \label{eqn:steiner_conceptual}
  \vol[d](A_\ge^{ext}) = \sum_{\i \in \{0,1,\dots,d-1\}} c_\i \ge^{d-\i}.
\end{align}
and it is clear that Corollary~\ref{cor:self-similar,_one_generator} may be rewritten
  \linenopar
\begin{align}
  \label{eqn:main_cor_conceptual}
  V_\tiling(\ge)
  &=  \sum_{\gw \in \Ds \cup \intDim} \negsp[1] c_\gw \ge^{d-\gw},
\end{align}
where for each fixed $\gw \in \Ds$,
  \linenopar
\begin{align}
  \label{def:c_wj}
  c_\gw
  &:=  \res{\gzs(s)} \sum_{\i=0}^d \frac{\genir^{\gw-\i}}{\gw-\i} \crv_\i.
\end{align}
Note that when $\gw = \i \in \{0,1,\dots,d-1\}$, one has $c_\gw = c_\i = \gzs(\i) \crv_\i$.
The obvious similarities between the tube formulas \eqref{eqn:steiner_conceptual} and \eqref{eqn:main_cor_conceptual} is striking. Our tube formula is a fractal power series in \ge, rather than just a polynomial in \ge (as in Steiner's formula). Moreover, our series is summed not just over the `integral dimensions' \intDim, but also over the countable set \Ds of scaling complex dimensions. The \coeffs $c_\gw$ of the tube formula are expressed in terms of the `curvatures' and the inradii of the generators of the tiling.

\begin{remark}\label{rem:obtaining-the-ext-from-the-int}
  The two formulas \eqref{eqn:steiner_conceptual} and \eqref{eqn:main_cor_conceptual} initially appear to be measuring very different things, but this may be misleading. It can happen that the exterior \epsnbd of the fractal itself is, in fact, equal to the union of the inner \epsnbd of the tiling and the exterior \epsnbd of its convex hull. In this case, \eqref{eqn:steiner_conceptual} and \eqref{eqn:main_cor_conceptual} are closely related. In fact, the tube formula for the exterior \ge-\nbd of the fractal is precisely the sum of the tube formula for the tiling and Steiner's tube formula for its convex hull. This occurs, for example, with the Sierpinski gasket; see Remark~\ref{rem:compatibility} and Figure~\ref{fig:gasket-compatibility}.
\end{remark}

\begin{remark}[Comparison of $V_\tiling$ with the Steiner formula]
  \label{rem:comparison with Steiner}
  In the trivial situation when the spray consists only of finitely many scaled copies of a diphase generator (so the scaling measure is supported on a finite set), the scaling zeta \fn will have no poles in \bC, so that $\Ds = \es$. Therefore, the tube formula becomes a sum over only the numbers $\i=0,1,\dots,d-1$,\footnote{It turns out during the computation that the \nth[d] summand always vanishes, so the sum extends only up to $d-1$.} for which the residues
  simplify greatly. In fact, in this case $\gzh(\i) = \gr_1^\i + \dots + \gr_N^\i$, so each residue from \eqref{eqn:single_gen_result} is a finite sum
  \linenopar
  \begin{align*}
    \gzh(\i)\crv_{\i}(\ge)
    &= \gr_1^\i \crv_{\i} \ge^{d-\i} + \dots + \gr_N^\i \crv_{\i} \ge^{d-\i}
    = \crv_{\i} (r_{w_1} \gen) \ge^{d-\i}
    + \dots + \crv_{\i} (r_{w_N} \gen) \ge^{d-\i}
  \end{align*}
  where $N$ is the number of scaled copies of the generator \gen, and $r_{w_n}$ is the corresponding scaling factor. Summing over \i, we obtain a tube formula for the scaled basic shape $r_{w_n} \gen$, for each $n=1,\dots,N$. The pluriphase case is analogous.
\end{remark}

\pgap

\section{Two illustrative examples}
\label{sec:examples}

\begin{exm}[The Sierpinski gasket tiling]
  \label{exm:Sierpinski gasket tiling}

  \begin{figure}
    \centering
    \includegraphics{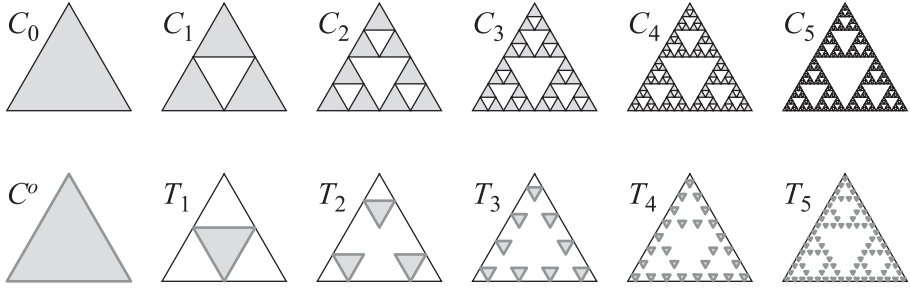}
    \caption{The Sierpinski gasket tiling. The first stages of the construction of the gasket are labeled $\hull_n$. The tiles $\{\simt_w(\gen)\}$ are labeled $\tileset_n$ where $|w|=n$.}
    \label{fig:Sierpinski gasket tiling}
  \end{figure}

  The Sierpinski gasket tiling \sSG (see Fig.~\ref{fig:Sierpinski
  gasket tiling}) is constructed via the system
  \linenopar
    \[\simt_1(z) := \tfrac12 z, \q
      \simt_2(z) := \tfrac12 z + \tfrac12, \q
      \simt_3(z) := \tfrac12 z + \tfrac{1+\ii\sqrt3}4,\]
  which has one common scaling ratio \(\simrat = 1/2,\) with $J=3$ and one generator
  $G$: an equilateral triangle of side length $\tfrac12$ and inradius
  \(\genir = \tfrac1{4\sqrt3}.\) Thus \sSG has inradii $\gr_k = \genir
  \simrat^{k}$ with multiplicity $3^k$, so the scaling zeta \fn is
  \linenopar
  \begin{align}
    \gzs(s) = \frac1{1 - 3 \cdot 2^{-s}},
  \end{align}
  and the scaling complex dimensions are
  \linenopar
  \begin{align}
    \Ds = \{D + \ii n\per \suth n \in \bZ\}
    \qq \text{for } D=\log_23, \; \per=\tfrac{2\gp}{\log2}.
  \end{align}
  The tube formula for \sSG is readily computed to be
  \linenopar
  \begin{align}
    V_{\sSG}(\ge)
    &= \tfrac{\sqrt3}{16 \log2} \sum_{n \in \bZ}
      \left(-\tfrac1{D+\ii n\per} + \tfrac{2}{D-1+\ii n\per} - \tfrac{1}{D-2+\ii n\per}\right)\left(\tfrac{\ge}{\genir}\right)^{2-D-\ii n\per}
      + \tfrac{3^{3/2}}{2} \ge^2 - 3 \ge.
    \notag 
  \end{align}
\end{exm}

\begin{exm}[The Koch tiling]
  \label{exm:The Koch Tiling}

  The standard Koch tiling \sK (see Figure~\ref{fig:koch-tiled-in}) is constructed via the self-similar system given in \eqref{eqn:koch-system}. The attractor of $\{\simt_1,\simt_2\}$ is the classical von Koch curve $K$, as in Figure~\ref{fig:koch-snowflake}. This system has one scaling ratio \(\simrat=|\gx| = 1/\sqrt3,\) with $J=2$ and one generator $G$: an equilateral triangle of side length $\tfrac13$ and generating inradius \(g=\tfrac{\sqrt3}{18}\). This tiling has inradii $\gr_k = g\simrat^{k}$ with multiplicity $2^k$, so the scaling zeta \fn is
  \linenopar
  \begin{align}
    \gzs(s) = \frac1{1 - 2 \cdot 3^{-s/2}},
  \end{align}
  and the scaling complex dimensions are
  \linenopar
  \begin{align}
    \Ds = \{D + \ii n\per \suth n \in \bZ\}
    \qq \text{for } D=\log_34, \; \per=\tfrac{4\gp}{\log3}.
  \end{align}
  By inspection, a tile with inradius $1/x$ will have tube formula
  \linenopar
  \begin{align}
    \label{eqn:Koch gtf}
    \gtf(x,\ge) =
      \begin{cases}
        3^{3/2} \left(-\ge^2 + 2\ge x\right), &\ge \leq 1/x, \\
        3^{3/2} x^2, &\ge \geq 1/x.
      \end{cases}
  \end{align}
  For fixed $x$, \eqref{eqn:Koch gtf} is clearly \cn at $\ge = 0^+$. Thus we have $\crv_{0} = -3^{3/2}, \crv_{1} = 2\cdot 3^{3/2}$, and $\crv_{2} = -3^{3/2}$. Now applying \eqref{eqn:single_gen_result}, the tube formula for the Koch tiling \sK is
  \linenopar
  \begin{align}
    V_{\sK}(\ge)
    &= \frac{\genir}{\log3} \sum_{n \in \bZ}
      \left(-\tfrac1{D+\ii n\per} + \tfrac{2}{D-1+\ii n\per} - \tfrac{1}{D-2+\ii n\per}\right)\left(\tfrac{\ge}{\genir}\right)^{2-D-\ii n\per} \notag \\
      &\hstr[12]+ 3^{3/2} \ge^2 + \tfrac1{1-2\cdot3^{-1/2}} \ge,
    \label{eqn:Koch final in eps/genir}
  \end{align}
  where $D=\log_34$, \(g=\tfrac{\sqrt3}{18}\) and
  $\per=\tfrac{4\gp}{\log3}$ as before.

  \begin{remark}[Nonlattice Koch tilings]
    \label{rem:nonlattice Koch}
    By replacing \(\gx = \tfrac{1}{2} + \tfrac{1}{2\sqrt{3}}\ii\) in \eqref{eqn:koch-system} with any other complex number satisfying
    \(|\gx|^2 + |1-\gx|^2 < 1,\)
    one obtains a family of examples of nonlattice self-similar tilings, as illustrated in Figure~\ref{fig:skew-koch}. Computation of the tube formula does not differ significantly from the lattice case. Further discussion of nonlattice Koch tilings may be found in \cite{SST}. As discussed in \cite{TFCD}, for example, this furnishes a family of tilings, almost all of which are Minkowski measurable.
    \begin{figure}
      \centering
      \includegraphics{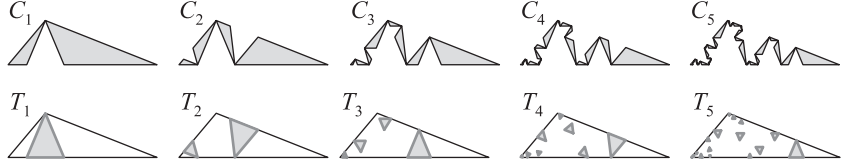}
      \caption{A measurable Koch tiling obtained by picking a nonlattice parameter \gx in \eqref{eqn:koch-system}.}
      \label{fig:skew-koch}
    \end{figure}
  \end{remark}

  \begin{remark}
    \label{rem:comparison of Koch tiling tube to Koch tube}
    In \S\ref{sec:koch} we discussed how to obtain a tube formula was obtained for the \ge-\nbd of the Koch curve itself (rather than of the tiling associated with it) and this led to the prediction that the complex dimensions of the curve are
  \linenopar
    \[\D_{\sK\star} = \{D+\ii n\per \suth n \in \bZ\} \cup \{0+\ii n\per \suth n \in \bZ\},\]
    where $D=\log_34$ and $\per=\frac{2\gp}{\log3}$. The line of poles above $D$ was expected\footnote{This set of complex dimensions was predicted in \cite{FGNT}, \S10.3, except for the dimensions above 0.}, and agrees precisely with the results of \S\ref{sec:tilings}. The meaning of the line of poles above $0$ is still unclear. We invite the reader to compare \eqref{eqn:Koch final in eps/genir} with the formula for the Koch curve. When formula \eqref{eqn:preview-of-corollary} for the area of the inner \ge-\nbd of the Koch snowflake is stated in full detail, it appears as follows.
  \linenopar
      \begin{equation}
        \label{eqn:final-result-compare}
        V(\ge) = G_1(\ge) \ge^{2-D} + G_2(\ge) \ge^2,
      \end{equation}
      where $G_1$ and $G_2$ are periodic \fns of multiplicative
      period 3, given by
  \linenopar
      \begin{subequations}
        \label{eqngrp:the periodic fns-compare}
        \begin{align}
          \label{eqn:first periodic fn-compare}
          G_1(\ge) :&= \frac1{\log3} \sum_{n \in \bZ}
            \left(a_n + \sum_{\ga \in \bZ} b_\ga g_{n-\ga} \right)
            (-1)^n \ge^{-in\per}
          \\
          \label{eqn:second periodic fn-compare}
          \text{ and \;} G_2(\ge) :&= \frac1{\log3} \sum_{n \in \bZ}
            \left(\gs_n + \sum_{\ga \in \bZ} \gt_\ga g_{n-\ga} \right)
            (-1)^n \ge^{-in\per},
        \end{align}
      \end{subequations}
    where $g_\ga$ are the Fourier \coeffs of the multiplicatively periodic \fn $h(\ge)$ discussed just before the statement of Theorem~\ref{thm:preview-of-final-result}, and
    $a_n, b_n, \gs_n,$ and $\gt_n$ are
    the complex numbers given by
  \linenopar
    \begin{align}
      a_n &= \frac{\gp - 3^{3/2}}{2^3(D + in\per)}
      +\frac{3^{3/2}}{2^3(D - 1 + in\per)}
      -\frac{3^{5/2}}{2^5(D - 2 + in\per)}
      +\frac12 b_n,
      \notag \\
      b_n &= \sum_{m=1}^\iy
       \frac{ (2m)! \; (3^{2m+1}-4)}
         {4^{2m+1} (m!)^2 (4m^2-1) (3^{2m+1}-2) (D - 2m - 1 + in\per)},
      \notag \\
      \label{eqn:mean-coeffs-compare}
      \gs_n &= - \log3\left(\frac{\gp}{3} + 2\sqrt{3}\right)\gd_0^n
        - \gt_n, \text{ and } \\
      \gt_n &= \sum_{m=1}^\iy
       \frac{(2m)! \; (3^{2m+1}-1)}
         {4^{2m-1} (m!)^2 (4m^2-1) (3^{2m+1}-2) (-2m - 1 + in\per)}.
      \notag
    \end{align}
    In the definition of $\gs_n$ we use the Kronecker delta to indicate a term that only appears in $\gs_0$.
  \end{remark}
\end{exm}

\begin{figure}
    \centering
    \includegraphics{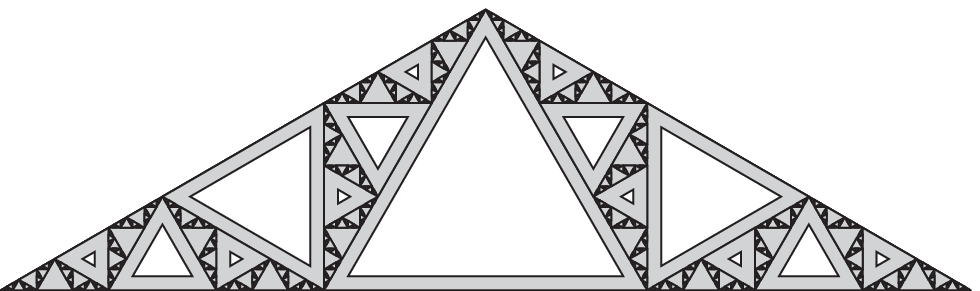}
    \caption{The exterior \ge-\nbd of the Koch curve is not simply related to the exterior \ge-\nbd of $\hull=[K]$ and the inner \ge-\nbd of the Koch tiling \sK.}
    \label{fig:koch-tile-inner-nbd}
  \end{figure}

\begin{remark}\label{rem:compatibility}
  As mentioned in Remark~\ref{rem:obtaining-the-ext-from-the-int}, the exterior \ge-\nbd of the Sierpinski gasket is the disjoint union of the inner \ge-\nbd of the tiling and the exterior \ge-\nbd of the largest triangle; see Figure~\ref{fig:gasket-compatibility}.
  \linenopar
  \begin{figure}
    \centering
    \includegraphics{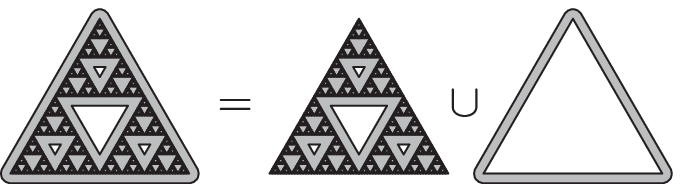}
    \caption{The exterior \ge-\nbd of the Sierpinski gasket is the disjoint union of the exterior \ge-\nbd of $\hull=[\attr]$ and the inner \ge-\nbd of the Sierpinski gasket tiling.}
    \label{fig:gasket-compatibility}
  \end{figure}
  This means that one can immediately obtain the tube formula for the exterior \ge-\nbd of the gasket by adding the tube formula for the tiling \sSG and Steiner's tube formula for its convex hull, the equilateral triangle $\hull_0$ (as labeled in Figure~\ref{fig:Sierpinski gasket tiling}).
  However, things do not always work out so neatly, as the example of the Koch tiling shows; see Figure~\ref{fig:koch-tile-inner-nbd}.
  In the forthcoming paper \cite{GeometryOfSST}, precise \conds are given for equality to hold as in \eqref{eqn:extnbd_decomp}. In particular,
  \linenopar
  \begin{align}
    \vol[d](\tiling_\ge^{ext}) = V_\tiling(\ge) + \vol[d]([\attr])
    \label{eqn:extnbd_decomp}
  \end{align}
  holds whenever one of the following equivalent conditions is verified:
  \begin{enumerate}
    \item $\del \tiling = \attr$.
    \item $\del \gen_q \ci \attr$ for all $q=1,\dots,Q$.
    \item $\del (\hull \less \simt(\hull)) \ci \attr$.
    \item $\del \hull \ci \attr$.
    \item $\attr_\ge^{ext} \cap \hull = T_\ge$ for some (and, equivalently, all) $\ge \geq 0$.
    \item $\attr_\ge^{ext} \cap \hull^c = \hull_\ge^{ext} \less \hull$ for some (and, equivalently, all) $\ge \geq 0$.
  \end{enumerate}
  Thus \eqref{eqn:main-result} (in conjunction with \eqref{eqn:extnbd_decomp}) allows one to compute explicit tube formulas for a large family of self-similar sets.
\end{remark}

\begin{exm}[The pentagasket tiling] 
  \label{exm:Pentagasket Tiling}

  The pentagasket tiling \sP is constructed via the self-similar system defined by the five maps
  \linenopar
  \[\simt_\j(x) = \tfrac{3-\sqrt5}2 x + p_\j, \qq \j=1,\dots,5,\]
  with common scaling ratio $\simrat=\phi^{-2}$, where $\phi = (1+\sqrt5)/2$ is the golden ratio, and the points $\frac{p_\j}{1-r} = c_\j$ form the vertices of a regular pentagon of side length 1; see Fig.~\ref{fig:pentagasket tiling}.
  \linenopar
  \begin{figure}
    \centering
    \includegraphics{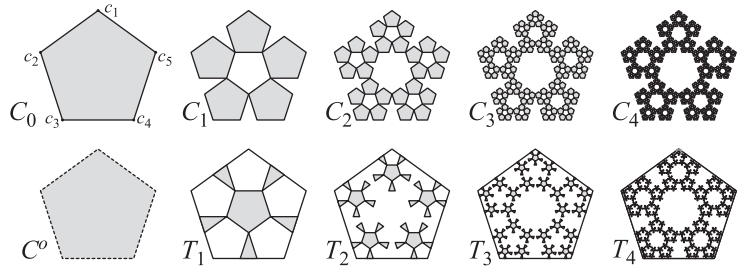}
    \caption{The pentagasket tiling \sP.}
    \label{fig:pentagasket tiling}
  \end{figure}

  The pentagasket \sP provides an example of multiple (noncongruent) generators $\gen_q$ with $q=1,2,\dots,6$. Specifically, $\gen_1$ is a regular pentagon and $\gen_2,\dots,\gen_6$ are congruent isosceles triangles, as seen in $\tileset_1$ of Fig.~\ref{fig:pentagasket tiling}. This example is developed fully in \cite[\S9.4]{TFCD}.
\end{exm}

\pgap

\bibliographystyle{harvard}

\end{document}